\documentstyle[12pt]{amsart}
\begin{document}
  \title{Harmonic morphisms of Allof-Wallach spaces of  positive curvature}
  \title[Harmonic morphisms of Allof-Wallach spaces]
  {Harmonic morphisms of Allof-Wallach spaces of positive curvature}
   \author{Hajime Urakawa}
  \address{Graduate School of Information Sciences, 
  Tohoku University, Aoba 4-3-09, Sendai 980-8579, Japan}
  \curraddr{}
  \email{urakawa@@math.is.tohoku.ac.jp}
  \thanks{
  Supported by the Grant-in-Aid for the Scientific Research, (C) No. 25400154, Japan Society for the Promotion of Science. 
  }
\maketitle
\begin{abstract} 
An infinite family of 
distinct harmonic morphisms  with minimal circle fibers
from the $7$-dimensional 
homogeneous Allof-Wallach spaces 
of positive curvature onto the $6$-dimensional flag manifolds is given.
    \end{abstract}
\numberwithin{equation}{section}
\theoremstyle{plain}
\newtheorem{df}{Definition}[section]
\newtheorem{th}[df]{Theorem}
\newtheorem{prop}[df]{Proposition}
\newtheorem{lem}[df]{Lemma}
\newtheorem{cor}[df]{Corollary}
\newtheorem{rem}[df]{Remark}
\section{Introduction}
In 1965, Eells and Sampson \cite{ES} initiated 
a theory of harmonic maps in which  
variational problems play central roles in geometry;\,Harmonic map is one of solutions of important variational problems which is a critical point of the energy functional 
$E(\varphi)=\frac12\int_M\vert d\varphi\vert^2\,v_g$ 
for smooth maps $\varphi$ of $(M,g)$ into $(N,h)$. The Euler-Lagrange equations are given by the vanishing of the tension filed 
$\tau(\varphi)$. 
On the other hand, 
Fuglede \cite{F} in 1978 and Ishihara \cite{Is} in 1979, introduced independently 
the alternative notion of harmonic morphism which 
preserves harmonic functions 
(see \cite{BW}). Harmonic morphisms are one of important examples of harmonic maps. 
\par
In this paper, we give new examples of harmonic morphisms, indeed, an infinite family of distinct harmonic morphisms from the $7$ dimensional homogeneous space of positive sectional curvature into the $6$ dimensional flag manifold. Namely, we show the following theorem.
\vskip0.6cm\par
\textbf{Theorem}  (cf. {\it Theorem 4.2}) 
{\em
Let $(P,g)=(M_{k,\,\ell},g_t)=(SU(3)/T_{k,\ell},g_t)$,  
$k,\,\,\ell\in {\mathbb Z},\,\,(k,\,\ell)=1;\,\,-1<t<0$, or $0<t<\frac13$, 
be infinitely many distinct homogeneous the $7$-dimensional Allof-Wallach spaces 
of positive sectional curvature,  and let $(M,h)$, 
the $6$-dimensional flag manifold $(SU(3)/T,h)$. 
Then, the Riemannian submersions with circle fibers,  
$\pi:\,(P,g)=(M_{k,\,\ell},g_t)\rightarrow (M,h)=(SU(3)/T,h)$,   are all harmonic morphisms with minimal fibers. }
\vskip0.6cm\par
Here,  the subgroups $T_{k,\,\ell}$ and $T$ of $SU(3)$ 
and the homogeneous space $M_{k,\,\ell}$ are given as follows. 
\vskip0.3cm\par
\begin{align*}
T_{k,\,\ell}&=\left\{
\begin{pmatrix}
e^{2\pi ki\theta}&0&0\\
0&e^{2\pi i \ell\theta}&0\\
0&0&e^{-2\pi i(k+\ell)}
\end{pmatrix}\big\vert\,\,\theta\in {\mathbb R}
\right\}
\\
&\subset 
T=
\left\{
\begin{pmatrix}
e^{2\pi i\theta_1}&0&0\\
0&e^{2\pi i\theta_2}&0\\
0&0&e^{-2\pi i(\theta_1+\ell_2)}
\end{pmatrix}\big\vert\,\,\theta_1, \,\,\theta_2\in {\mathbb R}
\right\}
\subset G=SU(3), \\
\mbox{and}\,\,&M_{k,\,\ell}=G/T_{k,\,\ell}=SU(3)/T_{k,\ell}.
\end{align*}  
\vskip0.6cm\par
\section{Preliminaries}
\subsection{Harmonic maps.}
We first prepare the materials for the first and second variational formulas for the bienergy functional and biharmonic maps. 
Let us recall the definition of a harmonic map $\varphi:\,(M,g)\rightarrow (N,h)$, of a compact Riemannian manifold $(M,g)$ into another Riemannian manifold $(N,h)$, 
which is an extremal 
of the {\em energy functional} defined by 
$$
E(\varphi)=\int_Me(\varphi)\,v_g, 
$$
where $e(\varphi):=\frac12\vert d\varphi\vert^2$ is called the energy density 
of $\varphi$.  
That is, for any variation $\{\varphi_t\}$ of $\varphi$ with 
$\varphi_0=\varphi$, 
\begin{equation}
\frac{d}{dt}\bigg\vert_{t=0}E(\varphi_t)=-\int_Mh(\tau(\varphi),V)v_g=0,
\end{equation}
where $V\in \Gamma(\varphi^{-1}TN)$ is a variation vector field along $\varphi$ which is given by 
$V(x)=\frac{d}{dt}\big\vert_{t=0}\varphi_t(x)\in T_{\varphi(x)}N$, 
$(x\in M)$, 
and  the {\em tension field} is given by 
$\tau(\varphi)
=\sum_{i=1}^mB(\varphi)(e_i,e_i)\in \Gamma(\varphi^{-1}TN)$, 
where 
$\{e_i\}_{i=1}^m$ is a locally defined orthonormal frame field on $(M,g)$, 
and $B(\varphi)$ is the second fundamental form of $\varphi$ 
defined by 
\begin{align}
B(\varphi)(X,Y)&=(\widetilde{\nabla}d\varphi)(X,Y)\nonumber\\
&=(\widetilde{\nabla}_Xd\varphi)(Y)\nonumber\\
&=\overline{\nabla}_X(d\varphi(Y))-d\varphi(\nabla_XY),
\end{align}
for all vector fields $X, Y\in {\frak X}(M)$. 
Here, 
$\nabla$, and
$\nabla^h$, 
 are Levi-Civita connections on $TM$, $TN$  of $(M,g)$, $(N,h)$, respectively, and 
$\overline{\nabla}$, and $\widetilde{\nabla}$ are the induced ones on $\varphi^{-1}TN$, and $T^{\ast}M\otimes \varphi^{-1}TN$, respectively. By (2.1), $\varphi$ is {\em harmonic} if and only if $\tau(\varphi)=0$. 
\vskip0.3cm\par
\subsection{Riemannian submersions.} 
We prepare with several notions on the Riemannian submersions. 
A $C^{\infty}$ mapping $\pi$ of a $C^{\infty}$ Riemannian manifold $(P,g)$ into another $C^{\infty}$ Riemannian manifold $(M,h)$ is called a {\em Riemannia submersion} if 
$(0)$ $\pi$ is surjective, $(1)$ the differential 
$d\pi=\pi_{\ast}:\,\,T_uP\rightarrow T_{\pi(u)}M\,\,(u\in P)$ 
of $\pi:\,\,P\rightarrow M$ 
is surjective for each $u\in P$, and 
$(2)$ each tangent space $T_uP$ at $u\in P$ has 
the direct decomposition: 
$$
T_uP={\mathcal V}_u\oplus {\mathcal H}_u,\qquad (u\in P), 
$$
which is orthogonal decomposition with respect to $g$ such  that ${\mathcal V}={\mbox{\rm Ker}}(\pi_{\ast\,u})\subset T_uP$ and 
$(3)$ the restriction of the differential $\pi_{\ast}=d\pi_u$ to 
${\mathcal H}_u$ is a surjective isometry,  
$\pi_{\ast}:\,\,({\mathcal H}_u,g_u)\rightarrow (T_{\pi(u)}M, h_{\pi(u)})$ for each $u\in P$ (cf. \cite{BW}). 
A manifold $P$ is the total space of a Riemannian submersion over $M$ with the  projection $\pi:\,\,P\rightarrow M$ onto $M$, 
where $p=\dim P=k+m$, $m=\dim M$, and $k=\dim\pi^{-1}(x)$, $(x\in M)$. 
A Riemannian metric $g$ on $P$, called 
{\em adapted metric } on $P$ which satisfies  
\begin{equation}
g=\pi^{\ast}h+k
\end{equation}
where $k$ is the Riemannian metric on each fiber 
$\pi^{-1}(x)$, $(x\in M)$. Then, 
$T_uP$ has the orthogonal direct decomposition 
of the tangent space $T_uP$, 
\begin{equation}
T_uP={\mathcal V}_u\oplus {\mathcal H}_u,\qquad\quad u\in P, 
\end{equation}
where the subspace ${\mathcal V}_u=\mbox{\rm Ker}(\pi_{\ast}{}_u)$ at $u\in P$, 
the {\em vertical subspace}, and the subspace 
${\mathcal H}_u$ of $P_u$ is called 
{\em horizontal subspace} at $u\in P$ which is the orthogonal complement of ${\mathcal V}_u$ in $T_uP$ with respect to $g$. 
\par 
In the following, we fix a locally defined orthonormal frame field, called {\em adapted local orthonormal frame field} to the projection $\pi:\,\,P\rightarrow M$, 
$\{e_i\}_{i=1}^p$ corresponding to $(2.9)$ 
in such a way that 
\par $\bullet$
$\{e_i\}_{i=1}^m$ is a locally defined orthonormal basis of the horizontal 
\par\quad  subspace ${\mathcal H}_u$ $(u\in P)$, 
and 
\par $\bullet$ 
$\{ e_i\}_{i=1}^k$ 
is a locally defined orthonormal basis 
of the vertical \par\quad 
subspace ${\mathcal V}_u$ $(u\in P)$. 
\medskip\par
Corresponding to the decomposition $(2.9)$, the tangent vectors $X_u$, and $Y_u$ in $T_uP$ 
 can be defined by 
\begin{align}
&X_u=X_u^{{\rm V}}+X_u^{{\rm H}},\quad 
Y_u=Y_u^{{\rm V}}+Y_u^{{\rm H}},\\
&X_u^{{\rm V}},\,\,
Y_u^{{\rm V}}
\in {\mathcal V}_u,\quad 
X_u^{{\rm H}},\,\,
Y_u^{{\rm H}}\in {\mathcal H}_u
\end{align}
for $u\in P$. 
\par
Then, there exist a unique decomposition 
such that 
$$
g(X_u,Y_u)=h(\pi_{\ast}X_u,\pi_{\ast}Y_u)
+k(X^{\rm V}_u,Y^{\rm V}_u),
\quad X_u,\,\,Y_u\in T_uP,\,\,u\in P. 
$$
\vskip0.3cm\par
\subsection{The reduction of the harmonic map equation}
Hereafter, we treat with the above problem more precisely in the case $\dim(\pi^{-1}(x))=1,\,\,(u\in P,\,\,\pi(u)=x)$. 
Let $\{e_1,\,\,e_1,\,\,\ldots,\,\,e_m\}$ be an adapted local orthonormal frame field being $e_{n+1}=e_m$, vertical. 
The frame fields 
$\{ e_i:\,\,i=1,2,\ldots,n\}$ are 
the basic orthonormal frame field 
on $(P,g)$ 
corresponds to an orthonormal 
frame field $\{\epsilon_1,\,\,\epsilon_2,\,\,\ldots,\,\,\epsilon_n\}$ on $(M,g)$.  
Here, a vector field $Z\in {\frak X}(P)$ is {\em basic} if 
$Z$ is horizontal and $\pi$-related to a vector field 
$X\in {\frak X}(M)$. 
\par
In this section, we determine the biharmonic equation 
precisely in the case that 
$p=m+1=\dim P$, $m=\dim M$, and $k=\dim \pi^{-1}(x)=1$ 
$(x\in M)$. 
Since $[V,Z]$ is a vertical field on $P$ if $Z$ is basic and $V$ is vertical (cf. \cite{ON}, p. 461).  Therefore, 
for each $i=1,\ldots,n$, 
$[e_i,e_{n+1}]$ is vertical, so we can write as follows. 
\begin{equation}
[e_i,e_{n+1}]=\kappa_i\,e_{n+1}, \quad
i=1,\,\,\ldots,\,\,n
\end{equation}
where $\kappa_i\in C^{\infty}(P)$ ($i=1,\ldots,n$). 
For two vector fields 
$X,\,\,Y$ on $M$, let $X^{\ast},\,\,Y^{\ast}$, be the 
horizontal vector fields on $P$. 
Then, 
$[X^{\ast},Y^{\ast}]$ is a vector field on $P$ which is 
$\pi$-related to a vector field $[X,Y]$ on $M$ (for instance, 
\cite{U2}, p. 143).  Thus, for $i,\,\,j=1,\ldots,n$, 
$[e_i,e_j]$ is $\pi$-related to 
$[\epsilon_i,\epsilon_j]$, and we may write as 
\begin{equation}
[e_i,e_j]=\sum_{k=1}^{n+1}D^k_{ij}\,e_k,
\end{equation}
where $D^k_{ij}\in C^{\infty}(P)\,\,
(1\leq i,\,j\leq n;\,1\leq k\leq n+1)$.
\subsection{The tension field} 
In this subsection, we calculate the tension 
field $\tau(\pi)$.  
We show that 
\begin{equation}
\tau(\pi)=-d\pi\left(\nabla_{e_{n+1}}e_{n+1}
\right)=-\sum_{i=1}^n\kappa_i\,\epsilon_i.
\end{equation}
Indeed, we have 
\begin{align}
\tau(\pi)&=\sum_{i=1}^m\left\{\nabla^{\pi}_{e_i}d\pi(e_i)-d\pi\left(\nabla_{e_i}e_i
\right)
\right\}\nonumber\\
&=\sum_{i=1}^n\left\{\nabla^{\pi}_{e_i}d\pi(e_i)-d\pi\left(\nabla_{e_i}e_i
\right)
\right\}
+\nabla^{\pi}_{e_{n+1}}d\pi(e_{n+1})-d\pi\left(\nabla_{e_{n+1}}e_{n+1}
\right)\nonumber\\
&=-d\pi\left(\nabla_{e_{n+1}}e_{n+1}
\right)\nonumber\\
&=-\sum_{i=1}^n\kappa_i\epsilon_i.\nonumber
\end{align}
Because, 
for $i,\,j=1,\ldots,n$, $d\pi(\nabla_{e_i}e_j)=\nabla^h_{\epsilon_i}\epsilon_j$, and 
$\nabla^{\pi}_{e_i}d\pi(e_i)=\nabla^h_{d\pi(e_i)}d\pi(e_i)=\nabla^h_{\epsilon_i}\epsilon_i$. Thus, we have 
\begin{equation}
\sum_{i=1}^n\left\{
\nabla^{\pi}_{e_i}d\pi(e_i)-d\pi\left(\nabla_{e_i}e_i
\right)
\right\}=0.
\end{equation}
Since $e_{n+1}=e_m$ is vertical, $d\pi(e_{n+1})=0$, so that 
$\nabla^{\pi}_{e_{n+1}}d\pi(e_{n+1})=0$. 
\par
Furthermore, we have, 
by definition of the Levi-Civita connection, we have, 
for $i=1,\ldots,n$, 
$$
2g(\nabla_{e_{n+1}e_{n+1}},e_i)=2g(e_{n+1},[e_i,e_{n+1}])=2\kappa_{i},
$$
and $2g(\nabla_{e_{n+1}}e_{n+1},e_{n+1})=0$. Therefore, we have 
$$
\nabla_{e_{n+1}}e_{n+1}=\sum_{i=1}^n\kappa_ie_i,
$$
and then, 
\begin{equation}
d\pi\left(\nabla_{e_{n+1}}e_{n+1}
\right)=\sum_{i=1}^n\kappa_i\epsilon_i.
\end{equation}
Thus, we obtain (2.9). \qed
\vskip0.6cm\par
Thus, we obtain the following theorem: 
\medskip\par
\begin{th} Let $\pi:\,\,(P,g)\rightarrow (M,h)$ be a Riemannian submersion over $(M,h)$. Then, 
\par
The tension field $\tau(\pi)$ of $\pi$ is given by 
\begin{align}
\tau(\pi)=
-d\pi\left(\nabla_{e_{n+1}}e_{n+1}
\right)
=
-\sum_{i=1}^n\kappa_i\epsilon_i,
\end{align}
where $\kappa_i\in C^{\infty}(P)$, $(i=1,\ldots,n)$. 
\end{th}
\medskip\par
\section{Harmonic morphisms}
 \vskip0.3cm\par
 \begin{df}
 (1) A smooth map $\pi:\,(P,g)\rightarrow (M,h)$ is 
  {\em harmonic} if 
 the tension field vanishes, $\tau(\pi)=0$, and\vskip0.1cm\par
 (2) \, $\pi:\,(P,g)\rightarrow (M,h)$ 
 is a {\em harmonic morphism} (cf. \cite{BW}, p. 106) if, for every 
 harmonic function $f:\,\,(M,h)\rightarrow {\mathbb R}$, 
 the composition 
 $f\,\circ \,\pi:\,\,(P,g)\rightarrow {\mathbb R}$ 
 is also harmonic. 
 \par
 (3) $\pi:\,(P,g)\rightarrow (M,h)$ 
 is {\em horizontally weakly conformal} (cf. \cite{BW}, p. 46) if, the differential $\pi_{{\ast}\,\,p}:\,\,T_pP\rightarrow T_{\pi(p)}M$ is surjective, and  
 $$(\pi^{\ast}h)(X,Y)=\Lambda(p)\,g(X,Y)
 \qquad (X,\,Y\in {\mathcal H}_p),$$ 
 for some non-zero number $\Lambda(p)\not=0\,(p\in P)$. 
 Here, ${\mathcal H}_p$ is the horizontal subspace 
 of $T_pP$ 
 for the Riemannian submersion 
 $\pi:\,(P,g)\rightarrow (M,h)$ 
 satisfying that 
\begin{equation*}
\left\{
\begin{aligned}
 T_pP&={\mathcal V}_p\oplus {\mathcal H}_p,
 \\
{\mathcal V}_p&={\rm Ker}(\pi_{\ast\,\,p}).
\end{aligned}
\right.
  \end{equation*}
 \end{df}
 
\vskip0.3cm\par\noindent
It is well known (cf. \cite{BW}, p. 108) that 
\vskip0.3cm
\begin{th} 
(Fuglede, 1978 and Ishihara, 1979) 
Let $\varphi:\,\,(P,g)\rightarrow (M,g)$ be a Riemannian submersion. 
Then, it is harmonic morphism if and only if 
the following two conditions hold: 
\par\quad 
(i)\,\,\, $\varphi:\,(P,g)\rightarrow (M,h)$ is harmonic and 
\par\quad 
(ii) \,$\varphi:\,(P,g)\rightarrow (M,h)$ is horizontally weakly conformal.  
\end{th} 
\vskip0.6cm\par
\begin{cor} {\rm (cf. \cite{BW}, p. 123)} 
A Riemannian submersion $\varphi:\,\,(P,g)\rightarrow (M,g)$ is a harmonic morphism if and only if $\varphi$ has minimal fibers.
\end{cor}
\vskip0.6cm\par
{\bf Example.} \quad 
Let $G$ be a compact Lie group, and 
$K\subset H\subset G$, closed subgroups of $G$. 
It is well known (\cite{BW}, p. 123)  that, the natural projection 
$\pi:\,\,(G/K,g)\rightarrow (G/H,h)$ is a harmonic Riemannian submersion with totally geodesic fibers. 
Notice that both the Riemannian metrics 
$g$ on $G/K$ and $h$ on $G/H$ are the invariant 
Riemannian metrics on $G/K$ and $G/H$ are induced from the $\mbox{Ad}(G)$-invariant inner product $\langle\,,\,\rangle$ on the Lie algebra $\frak g$ of $G$.  
The Allof-Wallach's Riemannian metrics 
of positive sectional curvature in our main theorem,  Theorem 4.2, are invariant Riemannian metrics, but are different from the ones induced from 
the $\mbox{Ad}(G)$-invariant inner product $\langle\,,\,\rangle$ on the Lie algebra $\frak g$ of $G$.  
\vskip0.6cm\par
\section{Statement of main theorem}
We first prepare the setting of the Allof-Wallach theorem \cite{AW}. Let $G=SU(3)$, 
and 
\begin{align*}
T_{k,\,\ell}&=\left\{
\begin{pmatrix}
e^{2\pi ki\theta}&0&0\\
0&e^{2\pi i \ell\theta}&0\\
0&0&e^{-2\pi i(k+\ell)}
\end{pmatrix}\big\vert\,\,\theta\in {\mathbb R}
\right\}
\\
&\subset 
T=
\left\{
\begin{pmatrix}
e^{2\pi i\theta_1}&0&0\\
0&e^{2\pi i\theta_2}&0\\
0&0&e^{-2\pi i(\theta_1+\theta_2)}
\end{pmatrix}\big\vert\,\,\theta_1\,\,\theta_2\in {\mathbb R}
\right\}\\
&\subset G_1=\left\{
\begin{pmatrix}
x&0\\
0&\mbox{det}(x^{-1})
\end{pmatrix}\vert\,x\in U(2)
\right\}
\subset G=SU(3), 
\end{align*}
and the Lie algebras of $G,\,T_{k,\,\ell},\,T,\,G_1$ 
by $\frak g$, ${\frak t}_{k,\,\ell}$, ${\frak t}$, ${\frak g}_1$, respectively. 
Let the Ad$(G)$-invariant inner product 
$\langle\,,\,\rangle_0$ by 
\begin{align*}
\langle X,\,Y\rangle_0&:=-{\mbox{Re}}({\mbox{Tr}}(XY)),\qquad X,\,Y\in {\frak g}, \\
{\frak m}={\frak{g}_1}^{\perp}&:=
\left\{\begin{pmatrix}
0&0&z_2\\
0&0&z_1\\
-\overline{z}_2&-\overline{z}_1&0
\end{pmatrix}\vert\,\,z_1,\,\,z_2\in {\mathbb C}
\right\}, \\
{\frak t}_{k,\,\ell}&:=\left\{
\begin{pmatrix}
2\pi ik\theta&0&0\\
0&2\pi \ell\theta&0\\
0&0&-2\pi i(k+\ell)\theta
\end{pmatrix}\vert\theta\in {\mathbb R}
\right\},\\
V_1:=&{{\frak t}_{k,\,\ell}}^{\perp}\cap {\frak g}_1,\qquad 
V_2:={{\frak g}_1}^{\perp}={\frak m}, 
\end{align*}
and let 
$$
{\frak g}={\frak su}(3)={\frak t}_{k,\,\ell}\oplus V_1\oplus V_2,
$$
the orthogonal direct decomposition of $\frak g$ with respect to the inner product $\langle\,,\,\rangle_0$. 
For $-1<t<\infty$, let the new inner product 
$\langle\,,\,\rangle_t$ 
by 
\begin{equation}
\langle x_1+x_2,y_1+y_2\rangle_t:=
(1+t)\langle x_1,\,y_1\rangle_0+\langle x_2,\,y_2\rangle_0,
\end{equation}
where
$x_i,\,\,y_i\in V_i\,\, (i=1,\,2)$, 
and let $g_t$, the corresponding 
$G$-invariant Riemannian metric on the homogeneous space 
$G/T_{k,\,\ell}$.  Then, 
\vskip0.6cm\par
\begin{th} (Allof and Wallach \cite{AW}) 
The homogeneous space $(G/T_{k,\,\ell},\,g_t)$ 
corresponding to $(4.1)$ with 
$(-1<t<0)$ or $(0<t<\frac13)$ have strictly positive sectional curvature.  
\end{th}
\vskip0.6cm\par
We state our main theorem as follows: 
\begin{th}
Let $\pi$ be the Riemannian submersion of 
$(SU(3)/T_{k,\,\ell},g_t)$ onto $(SU(3)/T,h)$, 
where $(SU(3)/T,h)$ is a flag manifold with the $SU(3)$-invariant Riemannian metric $h$ corresponding the inner product $\langle\,,\,\rangle_0$ on $\frak g$. 
Then, it is a harmonic morphism, i.e., 
for every harmonic function on a neighborhood V in 
$SU(3)/T$, the composition 
$f\circ\pi$ is harmonic on 
a neighborhood $\pi^{-1}(V)$ in $SU(3)/T_{k,\,\ell}$, 
and also it has minimal fibers.  
\end{th}
\vskip0.6cm\par
\section{Proof of main theorem}
Here, in this section, we give a proof of Theorem 4.2.
\par
We take a basis $\{X_0,X_1,X_2\}$ of $V_1$ and 
the one $\{X_3,X_4,X_5,X_6\}$ of $V_2$ as follows: 
\begin{align*}
X_0&=\frac{i}{\sqrt{5\,\Gamma}}\begin{pmatrix}
2\ell+k&0&0\\
0&2m+\ell&0\\
0&0&2k+m
\end{pmatrix},\\
X_1&=\frac{1}{\sqrt{2}}
\begin{pmatrix}
0&1&0\\
-1&0&0\\
0&0&0
\end{pmatrix},\quad 
X_2=\frac{1}{\sqrt{2}}\begin{pmatrix}
0&i&0\\
i&0&0\\
0&0&0
\end{pmatrix},\\
X_3&=\frac{1}{\sqrt{2}}
\begin{pmatrix}
0&0&1\\
0&0&0\\
-1&0&0
\end{pmatrix},\quad 
X_4=\frac{1}{\sqrt{2}}\begin{pmatrix}
0&0&i\\
0&0&0\\
i&0&0
\end{pmatrix},\\
X_5&=\frac{1}{\sqrt{2}}
\begin{pmatrix}
0&0&0\\
0&0&1\\
0&-1&0
\end{pmatrix},\quad 
X_6=\frac{1}{\sqrt{2}}\begin{pmatrix}
0&0&0\\
0&0&i\\
0&i&0
\end{pmatrix}. 
\end{align*}
Here $\Gamma:=k^2+\ell^2+k\ell$,  
$m:=-k-\ell$, and 
$\{X_0,X_1,X_2,X_3,X_4,X_5,X_6\}$ is an orthonormal 
basis of $V_1\oplus V_2$ with respect to 
$\langle\,,\,\rangle_0$, and 
${\frak g}={\frak su}(3)={\frak t}_{k,\,\ell}\oplus V_1\oplus V_2$.  
\par
Then, the basis of $V_1\oplus V_2$ 
$$
\left\{
\frac{1}{\sqrt{1+t}}\,X_0,\,\,\frac{1}{\sqrt{1+t}}\,X_1,\,\,\frac{1}{\sqrt{1+t}}\, X_2,\,\,X_3,\,\,X_4,\,\,X_5,\,\,X_6
\right\}
\qquad (5.1) 
$$
is orthonormal 
with respect to the inner product 
$\langle\,,\,\rangle_t$, $(-1<t<\infty)$ in $(4.1)$. 
We denote $M_{k,\,\ell}:=SU(3)/T_{k,\,\ell}$ with 
$k$ and $\ell\in{\mathbb Z}$ with $(k,\ell)=1$,  
and the corresponding local unit orthonormal vector fields on $M_{k,\,\ell}:=SU(3)/T_{k,\,\ell}$ by 
$$
\qquad\qquad\qquad\quad   \left\{e_0^t,\,\,e_1^t,\,\,e_2^t,\,\,e_3^t,\,\,e_4^t,\,\,e_5^t,\,\,e_6^t
\right\}\qquad\qquad\qquad\qquad   (5.2)
$$
\par
For the projection $\pi:\,M_{k,\,\ell}=SU(3)/T_{k,\,\ell}\rightarrow M=SU(3)/T$, 
each element $e^t_i$ $(i=0,\,1,\,\ldots,\,6)$ in (5.2) 
corresponds 
by $\pi_{\ast}$ (the differential 
of $\pi$),  
to each in 
$$
\left\{
0,\,\,\frac{1}{\sqrt{1+t}}\,e'_1,\,\,\frac{1}{\sqrt{1+t}}\, e'_2,\,\,e'_3,\,\,e'_4,\,\,e'_5,\,\,e'_6
\right\}, 
$$
where $\{e'_1,\,e'_2,\,e'_3,\,e'_4,\,e'_5,\,e'_6\}$ is an 
orthonormal frame field on $(SU(3)/T,h)$. 
\vskip0.6cm\par
By definition of the Levi-Civita connection of a Riemannian metric $g_t$, 
for every vector field $X$ on $P=M_{k,\,\ell}$, 
\begin{align*}
2g_t(X,\nabla^{g_t}_{e^t_0}e^t_0)&=
e^t_0\,g_t(X,e^t_0)+e^t_0\,g_t(X,e^t_0)-Xg_t(e^t_0,e^t_0)\\
&\qquad+g_t(e^t_0,[X,e^t_0])+g_t(e^t_0,[X,e^t_0])-g_t(X,[e^t_0,e^t_0])\\
&=2\left\{
e^t_0\,g_t(X,e^t_0)+g_t(e^t_0,[X,e^t_0])
\right\}.
\end{align*}
Thus, we have 
\begin{equation*}
\qquad 
g_t(X,\nabla^{g_t}_{e^t_0}e^t_0)=
e^t_0\,g_t(X,e^t_0)+g_t(e^t_0,[X,e^t_0]). 
\qquad (5.3) 
\end{equation*}
\par
On the other hand, 
we have 
$$
\qquad \qquad \quad 
e^t_0\,g_t(e^t_i,e^t_0)=0\qquad (i=0,1,\ldots,6)\qquad (5.4)
$$
Indeed, 
$$
\qquad \qquad 
e^t_0g(e^t_i,e^t_0)=0\qquad (i=0,1,\ldots,6)\quad \qquad (5.5)
$$
and by a straightforward computation, we have 
the following Lemma: 
\begin{lem}
We have 
\begin{align*}
&\left[
\frac{1}{\sqrt{t+1}}\,X_1,\,\frac{1}{\sqrt{t+1}}\,X_0
\right]=-\frac{3(k+\ell)}{(1+t)\sqrt{5\Gamma}}\,X_2,\\
&\left[
\frac{1}{\sqrt{t+1}}\,X_2,\,\frac{1}{\sqrt{t+1}}\,X_0
\right]
=\frac{3(k+\ell)}{(1+t)\sqrt{5\Gamma}}\,X_1,\\
&\left[X_3,\,\frac{1}{\sqrt{1+t}}\,X_0\right]
=\frac{-3\ell}{\sqrt{1+t}\,\sqrt{5\Gamma}}\,X_4,
\left[X_4,\,\frac{1}{\sqrt{1+t}}\,X_0\right]
=\frac{3\ell}{\sqrt{1+t}\,\sqrt{5\Gamma}}\,X_3,\\
&\left[X_5,\,\frac{1}{\sqrt{1+t}}\,X_0\right]
=\frac{3k}{\sqrt{1+t}\,\sqrt{5\Gamma}}\,X_6,
\left[X_6,\,\frac{1}{\sqrt{1+t}}\,X_0\right]
=\frac{-3k}{\sqrt{1+t}\,\sqrt{5\Gamma}}\,X_5.
\end{align*}
\end{lem}
\vskip0.8cm\par
By Lemma 5.1, we have 
\begin{align*}
g_t(e^t_0,[X,e^t_0])=0
\qquad (\forall\,\,X=X_i\,\,(i=0,1,\cdots,6)).\qquad (5.6)
\end{align*}
By (5.3),\,\,(5.4),\,\,(5.5), we have 
\begin{align*}
g_t(X,\nabla^{g_t}_{e^t_0}e^t_0)=0\qquad
(\forall X\in {\mathfrak X}(M_{k,\,\ell})), \qquad (5.7) 
\end{align*}
which implies that 
\begin{equation*}
\qquad\qquad\qquad\qquad \nabla^{g_t}_{e^t_0}e^t_0=0.\qquad\qquad\qquad\qquad\qquad  (5.8) 
\end{equation*}
Then, we have 
\begin{equation*}
\qquad\qquad\qquad \tau(\pi)=-d\pi(\nabla^{g_t}_{e^t_0}e^t_0
)=0.\qquad\qquad\qquad  (5.9)
\end{equation*}
Therefore, by (5.8) and (5.9), 
the submersion 
$\pi$ is a harmonic map with minimal fibers. 
Due to Corollary 3.3, we have Theorem 4.2. 
\qed
\vskip0.6cm\par
\begin{rem}
(1) In our main theorem, Theorem 4.2, since 
$(G/T,h)$ is a flag manifold, so a K\"{a}hler manifold, 
it admits a lot of harmonic function on an open 
subset $V$ in $G/T$. For a harmonic function 
$f$ on an open subset $V\subset G/T$, then 
$f\,\circ\,\pi$ is harmonic on $\pi^{-1}(V)$. 
 \par
 (2) Our fibration $\pi:\,\,(SU(3)/T_{k,\,\,\ell},g)\rightarrow (SU(3)/T,h)$ has close similarities to the Hopf fibration 
 $\pi':\,\,(S^{2n+1},g_0)\rightarrow ({\mathbb C}P^n,h_0)$. 
 Both the total spaces have positive sectional curvature, and both the base spaces are K\"{a}hler manifolds. 
\end{rem}
\vskip0.6cm\par

\end{document}